\documentclass[11pt]{article}

\newfont{\bb}{msbm10 at 11pt}
\def\r{\hbox{\bb R}}
\newcommand{\T}{\hbox{\bf t}}
\newcommand{\N}{\hbox{\bf n}}
\newcommand{\B}{\hbox{\bf b}}
\newcommand{\C}{\hbox{\bf c}}
\newcommand{\X}{\hbox{\bf X}}

\newtheorem{theorem}{Theorem}[section]
\newtheorem{corollary}[theorem]{Corollary}
\newtheorem{definition}[theorem]{Definition}
\newtheorem{remark}[theorem]{Remark}

\begin{document}

\title{On linear Weingarten  surfaces\footnote{Partially
supported by MEC-FEDER 
 grant no. MTM2004-00109.}}
\author{Rafael L\'opez\\
Departmento de Geometría y Topología\\
Universidad de Granada\\
18071 Granada (Spain)\\
e-mail:{\tt rcamino@ugr.es}}
\date{}
\maketitle


\begin{abstract} In this paper we study surfaces 
in Euclidean 3-space that satisfy a Weingarten condition of linear type as 
$\kappa_1=m \kappa_2 +n$, where $m$ and $n$ are real numbers
 and  $\kappa_1$ and $\kappa_2$ denote 
the principal curvatures at each point of the surface. We investigate the possible 
existence of such surfaces parametrized by a uniparametric family of circles. Besides 
 the surfaces of revolution, we prove that not exist more  except the 
case $(m,n)=(-1,0)$, that is, if the surface is one of the classical 
examples of minimal surfaces discovered by Riemann.
\end{abstract}

\section{Introduction}\label{intro}

A surface $S$ in  Euclidean 3-space $\r^3$ is called a {\it 
Weingarten surface} if there is some relation between its
two principal curvatures $\kappa_1$ and $\kappa_2$, 
that is to say, there is a smooth function $W$ of two variables such that
$W(\kappa_1,\kappa_2)=0$. 
The classification of the Weingarten surfaces in Euclidean space 
is almost completely open today. 
These  surfaces were introduced by the very Weingarten  \cite{we1,we2} 
in the context of the problem of finding all surfaces isometric to a given surface of revolution. Along the history they  have been of interest for geometers: for example,  
\cite{ch,hw,ho,vo} and more recently, \cite{gmm,ks,rs}. Applications of Weingarten surfaces  on computer aided design and  shape investigation  can seen in \cite{bg}.
In this work we study Weingarten surfaces that satisfy the simplest case for $W$, that is, that $W$ is  linear:
\begin{equation}\label{w1}
\kappa_1=m\ \kappa_2+ n
\end{equation}
where $m$ and $n$ are constants. We say then that $S$ is a {\it linear  Weingarten surface} and we abbreviate by a LW-surface. In particular, 
umbilical surfaces  $(m,n)=(1,0)$ or constant mean 
curvature surfaces ($m=-1$) are LW-surfaces. Throughout this work, we exclude the case that one of the principal curvatures is zero, that is, 
 we shall assume that $m\not=0$.

Among all LW-surfaces, the class of surfaces of revolution are particularly studied
 because
 in such case, 
Equation (\ref{w1}) leads to an ordinary differential equation. Its study is then simplified to find the profile curve that defines the surface \cite{ho}.  
On the other hand,  if  $S$ is a closed LW-surface of genus zero, it must be  a surface of revolution \cite{vo}. See generalizations in \cite{ks}.

The aim of this paper is the search of new LW-surfaces that generalize  the surfaces of revolution.   In this sense, we give the following 

\begin{definition} A cyclic surface in Euclidean space $\r^3$ is a surface determined by a smooth uniparametric family of pieces of circles.
\end{definition}

In particular, surfaces of revolution and tubes are cyclic surfaces. The motivation of the 
present work comes from what happens for the family of surfaces with constant mean curvature. When the mean curvature vanishes on the surface, that is, if the surface is minimal, $(m,n)=(-1,0)$, the only rotational minimal surface is the catenoid. Riemann found all non-rotational minimal surfaces foliated by circles in parallel planes \cite{ri}. 
Enneper proved that for a cyclic minimal surface, the planes containing the circles must be parallel \cite{en1,en2} and then, it is one of the examples obtained by Riemann. The Riemann examples  play a major  role in the theory of minimal surfaces. When the mean curvature is a non-zero constant, $(m,n)=(-1,n)$ and $n\not=0$, Nitsche proved that a cyclic surface must be a surface of revolution \cite{ni}, whose classification is well known \cite{de}.

In this paper, we study cyclic LW-surfaces. We call that a cyclic surface is of Riemann-type if the planes containing the circles 
of  the foliation are parallel. Our interest in this work is twofold. First, we want to know if  a cyclic LW-surface must be of Riemann-type. In this sense, we prove, 

\begin{quote} {\it S.1: A cyclic LW-surface with $(m,n)=(m,0)$ must be  of Riemann-type.}
\end{quote}

The restriction $n=0$ is merely technical since, as we will see, the proof involves long computations that in the case $n\not=0$ become very difficult to manage. However, we hope that the same result holds for the general case of $n\not=0$. 
On the other hand, and assuming now that the planes are parallel, we look for  new LW-surfaces. However, we conclude then that 

\begin{quote}{\it S.2: Besides the surfaces of revolution, the only LW-surfaces of 
Riemann-type with arbitrary pair $(m,n)$ are the classical Riemann examples of minimal surfaces, that is, if 
$(m,n)=(-1,0)$.}
\end{quote}

This gives a particularity of the Riemann examples in the family of 
LW-surfaces.

\begin{remark} Whenever we talk  of a LW-surface, we exclude the umbilical case, that is, 
$(m,n)=(-1,0)$. Moreover, we point out that any uniparametric family of (non-necessary parallel) planes intersects a sphere into circles. 

\end{remark}

\section{Preliminaries}

In this section, we fix some notation on local classical differential geometry of 
surfaces. Let $S$ be a surface in $\r^3$ and consider  $\X=\X(u,v)$ 
 a local parametrization of $S$ defined in the $(u,v)$-domain. 
Let $N$  denote the unit normal vector field on $S$ given by
$$N=\frac{\X_u\wedge \X_{v}}{|\X_{u}\wedge \X_{v}|},
\qquad \X_u=\frac{\partial \X}
{\partial u},\   \X_v=\frac{\partial \X}{\partial v},$$
where $\wedge$ stands the cross product of $\r^3$. 
In each tangent plane, the induced metric 
$\langle,\rangle$  is determined by the first
fundamental form
$$I=\langle d\X,d\X\rangle=Edu^2+2Fdudv+Gdv^2,$$
with differentiable coefficients
$$E=\langle \X_u,\X_u\rangle,\quad F=\langle\X_u,\X_v\rangle,\quad
G=\langle \X_v,\X_v\rangle.$$
The shape operator of the immersion is represented by the second fundamental
form
$$II=-\langle d N ,d\X\rangle=e\ du^2+2f\ dudv+g\ dv^2,$$
with 
$$e=\langle  N ,\X_{uu}\rangle,\quad f=\langle N ,\X_{uv}\rangle,
\quad g=\langle  N ,\X_{vv}\rangle.$$
Under this parametrization $\X$, the mean curvature $H$ and the Gauss curvature $K$ have the classical expressions
$$H=\frac{eG-2f F+gE}{2 (EG-F^2)},\hspace*{1cm}K=\frac{eg-f^2}{EG-F^2}.$$
Let us denote by $[, ,]$ the determinant in $\r^3$ and put
$W=EG-F^2$. Then $H$ and $K$ write  as  
\begin{equation}\label{mean}
H=\frac{G[\X_u,\X_v,\X_{uu}]-2F [\X_u,\X_v,\X_{uv}]+E [\X_u,\X_v,\X_{vv}]}{2W^{3/2}}:=\frac{H_1}{2W^{3/2}},
\end{equation}
\begin{equation}\label{gauss}
 K=\frac{[\X_u,\X_v,\X_{uu}] [\X_u,\X_v,\X_{vv}] -[\X_u,\X_v,\X_{uv}]^2}{W^2}:=\frac{K_1}{W^2}.
\end{equation}
The principal curvatures $\kappa_1$ and $\kappa_2$ are given then by 
$$\kappa_1=H+\sqrt{H^2-K},\hspace*{1cm}\kappa_2=H-\sqrt{H^2-K}.$$
Then the condition (\ref{w1}) writes now as
\begin{equation}\label{w2}
(1-m) H_1-2 W^{3/2} n=-(1+m)\sqrt{H_1^2-4W K_1}.
\end{equation}
After some manipulations, and squaring twice (\ref{w2}), we obtain
\begin{equation}\label{w3}
\left(-m H_1^2+(1+m)^2 W K_1+ n^2 W^3\right)^2-n^2(1-m)^2 H_1^2 W^3=0. 
\end{equation}

\section{Cyclic LW-surfaces}

In this section we prove the first statement $S.1$ of the Introduction, that is,

\begin{theorem} Let $S$ be a cyclic LW-surface with $(m,n)=(m,0)$. Then 
the planes of the foliation are parallel.
\end{theorem}

The methods that we apply in our proofs are based on \cite{ni}. Let $\Gamma=\Gamma(u)$ be an orthogonal smooth curve to each $u$-plane of the foliation and denote by $u$ its arc-length parameter. We assume that the planes of the foliation are not parallel and we shall arrive to a contradiction. 
Let $\T$ be the unit tangent vector to $\Gamma$. Consider the Frenet frame
of the curve $\Gamma$, $\{\T,\N,\B\}$, where $\N$ and $\B$ denote the normal and binormal vectors respectively. Locally we parametrize $M$ by
$$\X(u,v)=\C(u)+r(u)(\cos{v}\ \N (u)+\sin{v}\ \B (u)),$$
where $r=r(u)>0$ and $\C=\C(u)$ denote the radius and centre of each $u$-circle of the foliation. Consider the Frenet equations of the curve $\Gamma$:
\begin{eqnarray*}
\T'&=&\hspace*{.8cm}\kappa \N\\
\N'&=&-\kappa \T+\ \ \sigma \B\\
\B'&=&\hspace*{.6cm}-\sigma\N
\end{eqnarray*}
where the prime $'$ denotes the derivative with respect to the $u$-parameter and $\kappa$ and $\sigma$ are the curvature and torsion of $\Gamma$, respectively. Observe that $\kappa\not=0$ because $\Gamma$ is not a straight-line.
Also, set
\begin{equation}\label{alfa}
\C'=\alpha\T+\beta\N+\gamma\B,
\end{equation}
where $\alpha,\beta,\gamma$ are smooth functions on $u$.

By using the Frenet equations and (\ref{alfa}), a straightforward computation  shows that 
(\ref{w3}) can be expressed by trigonometric polynomial on $\cos{(jv)}$, $\sin{(jv)}$. Exactly, there exist smooth functions on $u$, namely $A_j$ and $B_j$, such that (\ref{w3}) writes as
\begin{equation}\label{formula}
A_0+\sum_{j=1}^6 (A_j \cos{(jv)} +B_j\sin{(j v)})=0.
\end{equation}
Since this is an expression on the independent trigonometric terms $\cos{(jv)}$ and $\sin{(jv)}$, all coefficients $A_i,B_i$ must vanish.
The values for $A_6$ and $B_6$ are:
$$A_6=-\frac{1}{32}(m-1)^2 \kappa^2 r^6\left(\beta^4+(\gamma^2-\kappa^2 r^2)^2+\beta^2(2\kappa^2 r^2-6\gamma^2)\right)=0.$$
$$B_6=-\frac{1}{8}(m-1)^2\beta\gamma\kappa^2 r^6(\beta^2-\gamma^2+\kappa^2 r^2)=0.$$
Recall that in the next reasoning, $\kappa\not=0$. 
From $B_6$ we consider three possibilities.

\begin{enumerate}
\item Case $\beta\gamma\not=0$. Then $\beta^2=\gamma^2-\kappa^2 r^2$. From $A_6=0$, we
obtain $-4\gamma^2(\gamma^2-\kappa^2 r^2)=0$. Since $\gamma\not=0$, then $\gamma^2=\kappa^2 r^2$. But then, $\beta=0$. As conclusion, this case is impossible.

\item Case $\gamma=0$. Then 
$$A_6=-\frac{1}{32}(m-1)^2\kappa^2 r^6(\beta^2+\kappa^2 r^2)^2,$$
which yields a direct contradiction.

\item Case $\beta=0$. Now
$$A_6=-\frac{1}{32}(m-1)^2\kappa^2 r^6(\gamma^2-\kappa^2 r^2)^2.$$
Hence $\gamma^2=\kappa^2 r^2$. Then 
$$A_4=-\frac18(6+m(6m-13))\kappa^4 r^8(\alpha^2-r'^2)=0.$$
$$B_4=\frac14(6+m(6m-13))\alpha\kappa^4 r^8 r'=0.$$

\begin{enumerate}
\item If $(6+m(6m-13))\not=0$, then $\alpha^2=r'^2$ and $\alpha r'=0$. Thus 
$\alpha=0$ and $r$ is a constant function. Then 
$A_2=-\frac12(2m^2-5m+2)r^{10}\kappa^6=0$. Then  $m=1/2$ and $m=2$. In both cases, the computations of $A_1$ gives $\tau=0$ and then (\ref{w3}) 
implies $\frac94 r^{10}\kappa^6$ and $9 r^{10}\kappa^6$ respectively. Anyway, we conclude a contradiction.

\item Therefore, it suffices to study the case that $(6+m(6m-13))=0$, that is, 
$m=2/3$ and $m=3/2$. For simplicity, we do the proof in the former case (the 
case $m=3/2$ is obtained interchanging the roles of $\kappa_1$ and 
$\kappa_2$ in the linear relation $\kappa_1=m\kappa_2$).

Before to follow, we point out that  the case $\alpha=0$ is impossible, because  
$$A_3=-\frac{5}{18}\kappa^3 r^8 r'^2\tau=0.$$ 
Then this means that  $r'=0$ or $\tau=0$. If $r$ is a constant function, $A_2=2/9 r^{10}\kappa^6=0$, which it is false. If $\tau=0$, then $B_2$ and $B_1$ give respectively, 
$$\kappa^3 r^2+2 r r'\kappa'-\kappa(9 r'^2+2 r r'')=0.$$
$$81\kappa^3 r^2+2 r r'\kappa'+\kappa(71 r'^2-2 r r'')=0.$$
By combining both equations, we have $80\kappa(\kappa^2 r^2+r'^2)=0$, which it is a contradiction.

From now, we assume $\alpha\not=0$.
Then the computation of $A_3$ and $B_3$ imply:
\begin{eqnarray*}x_1:&=&3\alpha^3\kappa-2\alpha^2\kappa r \tau-2\kappa r r'(\alpha'-\tau r')\\
&+&\alpha (\kappa^3 r^2+4 r \kappa' r'-\kappa(21 r'^2+2 r r''))=0.
\end{eqnarray*}
\begin{eqnarray*}
x_2&:=&\alpha^2(-2 r \kappa'+15\kappa r')+2\alpha \kappa r (\alpha'-2\tau r')\\
&+& r'(\kappa^3 r^2+2 r\kappa' r'-\kappa (9 r'^2+2r r''))=0.
\end{eqnarray*}
Then $\alpha x_1+r' x_2=0$ yields $(\alpha^2+r'^2)^2 x_3=0$, where 
$$x_3:=3\alpha^2\kappa+\kappa^3 r^2-2\alpha\kappa r\tau+2 r \kappa' r'-\kappa (9 r'^2 +2 r r'').$$
Since $\alpha\not=0$, then 
$x_3=0$. Now, $r' x_3-x_2=0$ implies 
$$x_4:=\alpha r \kappa'+\kappa (-6\alpha r'+r (-\alpha'+\tau r'))=0.$$
In this expression, we obtain $\kappa'$: 
$$\kappa'=\frac{\kappa (r\alpha'+6\alpha r'-r\tau r')}{\alpha r},$$
and substituting  into the value of $x_1$, we get
$$x_5:=3\alpha^3-2\alpha^2 r\tau+2 r r'(\alpha'-\tau r')+\alpha(\kappa^2 r^2+3 r'^2-2r r'')=0.$$
Hence we obtain the value of $r''$, which putting it into $A_2=0$ and $B_2=0$ give respectively 
\begin{eqnarray*}
y_1&:=& -7\alpha^4-8\alpha^2\kappa^2 r^2-\alpha^3 r\tau+\alpha\kappa^2 r^3\tau 
+78 \alpha^2 r'^2+10\kappa^2 r^2 r'^2\\
&-&\alpha r\tau r'^2-15 r'^4=0.
\end{eqnarray*}
$$y_2:= -42\alpha^3-18\alpha\kappa^2 r^2-\alpha^2 r \tau+\kappa^2 r^3 \tau+58\alpha r'^2-r \tau r'^2=0.$$
Now $y_1-\alpha y_2=0$ gives
$$7\alpha^2+2\kappa^2 r^2-3 r'^2=0.$$
From this equation, we obtain $r'^2$,
\begin{equation}\label{r2}
r'^2=\frac13\left(7\alpha^2+2\kappa^2 r^2\right)
\end{equation}
 and we introduce it into  $y_1=0$ 
concluding
$$280 \alpha^3+62\alpha\kappa^2 r^2-10\alpha^2 r\tau+\kappa^2 r^3\tau=0.$$
If $\kappa^2 r^2\not =10\alpha^2$, then 
$$\tau=-\frac{2(140\alpha^3+31\alpha\kappa^2 r^2)}{r(\kappa^2 r^2-10 \alpha^2)}.$$
By substituting in $A_1$ and using (\ref{r2}), we 
have 
$$45\alpha^4+31\alpha^2\kappa^2 r^2+5\kappa^4 r^4=0,$$
which it is a contradiction. 
As conclusion,  $\kappa^2 r^2 =10\alpha^2$ and $A_3=0$ gives
$$\tau=\frac{3\kappa^2 r^2+10 r'^2}{20 r r'}.$$
From (\ref{r2}), 
 we have $r'^2=9\alpha^2=9/10 \kappa^2 r^2$. Returning with the computations, 
the coefficient $A_2$ (or $B_2$) gives $\kappa^6 r^{10}=0$, obtaining  the desired contradiction.

\end{enumerate}
\end{enumerate}

\section{LW-surfaces of Riemann-type}

We consider a cyclic surface $S$ of Riemann type, that is, a cyclic surface
where the pieces of circles of the foliation lie in parallel planes, for example, 
parallel to the $x_1 x_2$-plane. Because our reasoning is local, we 
can assume that $S$ writes as 
$$X(u,v)=(a(u),b(u),u)+r(u)(\cos{v},\sin{v},0),$$
where  $a, b$ and $r$ are smooth function in some $u$-interval $I$ and $r>0$ denotes the 
radius of each circle of the foliation. Moreover, $S$ is a surface of 
revolution if and only if $a$ y $b$ are constant functions. 
If we compute (\ref{w3}), we obtain an expression 
\begin{equation}\label{p0}
\sum_{j=0}^{12} A_j(u) \cos{(j v)}+B_j(u)\sin{(j v)}=0.
\end{equation}
Again, the functions $A_j$ and $B_j$ on $u$ vanish on $I$. 
We distinguish two cases according to the value of $n$.

\begin{enumerate}
\item Case $n\not=0$.

The computation of   $A_{12}$ and $B_{12}$ give respectively:
$$A_{12}=\frac{1}{2048} n^4 r^{12} A\hspace*{1cm}
B_{12}=\frac{512}{n^4} r^{12}  B,$$
where
$$A=a'^{12} -66 a'^{10} b'^2+495 a'^8 b'^4-924 a'^6 b'^6+495 a'^4 b'^8- 
66 a'^2 b'^{10}+b'^{12}.$$
$$B=a'b'\left(3 a'^{10}-55 a'^8 b'^2+198 a'^6 b'^4-198 a'^4 b'^6+ 
55 a'^2 b'^8-3b'^{10}\right).$$
We assume now that $S$ is not a surface of revolution and we will arrive to a contradiction.
As $r>0$, $A=B=0$. Because the expressions of $A$ and $B$ do not depend on $r$, we do a 
change of variables. Since the planar curve $\alpha(u)=(a(u),b(u))$ is not constant, we reparametrize 
it by the length-arc, that is, $(a(u),b(u))=(x(\phi(u),y(\phi(u))$, where
\begin{equation}\label{para}
a'(u)=\phi'(u)\cos{(\phi(u))},\hspace*{.5cm}b'(u)=\phi'(u)\sin{(\phi(u))},
\hspace*{.5cm}\phi'^2=a'^2+b'^2.
\end{equation}
With this change, $A$ and $B$ write now as:
$$A=\phi'(u)^{12}\cos{(12\phi(u))},\hspace*{1cm}B=\phi'(u)^{12}\sin{(12\phi(u))}.$$
Therefore, $\phi'=0$, that is, $\alpha$ is   a constant curve: contradiction.

\item Case $n=0$.

Now (\ref{w3}) is simply $-m H_1^2+(1+m)^2 W K_1=0$ and
Equation (\ref{p0}) is then  a sum until $j=3$, with 
$$A_3=-\frac14(1+m)^2 r^5\left(a''(a'^2-b'^2)-2a'b'b''\right).$$
$$B_3=-\frac14(1+m)^2 r^5\left( b''(a'^2-b'^2)+2a'b'a''\right).$$
We assume that $S$ is not a surface with constant mean curvature, that is, $m\not=-1$. 
As in the case $n\not=0$, we assume that $S$ is not a surface of revolution and we will 
obtain a contradiction. As above,  we reparametrize the 
curve $\alpha(u)=(a(u),b(u))$ as in (\ref{para}).
Then $A_3=B_3=0$ lead to respectively:
$$\phi'(u)^2\left(-\phi''(u)\cos{(3\phi(u))}+\phi'(u)^2\sin{(3\phi(u))}\right)=0.$$
$$\phi'(u)^2\left(\phi'(u)^2\cos{(3\phi(u))}+\phi''(u)\sin{(3\phi(u))}\right)=0.$$
By combining both equations, we obtain $\phi'(u)=0$ on $I$, obtaining the desired contradiction.

\end{enumerate}

 In the case $m=-1$ and $n=0$,  $S$ is a  minimal surface. Then the degree of (\ref{p0}) is $2$. Here, $A_2=B_2=0$ imply 
\begin{equation}\label{ri1}
a'=\lambda r^2\hspace*{1cm} b'=\mu r^2
\end{equation}
for some  constants $\lambda,\mu\geq 0$.   Hence that (\ref{p0})  gives 
\begin{equation}\label{ri2}
1+(\lambda^2+\mu^2)r^4+r'^2-r r''=0.
\end{equation}
Equations (\ref{ri1}) and (\ref{ri2}) define the Riemann examples ($\lambda^2+\mu^2\not=0$) and
 the catenoid ($\lambda^2+\mu^2=0$).

As conclusion,

\begin{theorem} \label{t2} The only  LW-surfaces of Riemann-type are:
\begin{enumerate}
\item The surfaces of revolution.
\item The classical Riemann examples of minimal surfaces.
\end{enumerate}
\end{theorem}

Moreover, the Riemann examples can be viewed as an exceptional 
case in the family of cyclic LW-surfaces, at least with $n=0$, that is,

\begin{corollary} \label{co1} 
Riemann examples of minimal surfaces  are the only non-rotational cyclic 
surfaces that satisfy a linear Weingarten relation of type $\kappa_1=m\kappa_2$, $m\not=0$.
\end{corollary}

Theorem \ref{t2} and Corollary \ref{co1} show the statement $S.2$ of the Introduction.


\begin{thebibliography}{99}


\bibitem{bg}  van-Brunt B, K. Grant K (1996)
Potential applications of Weingarten surfaces in CAGD. I: Weingarten surfaces and surface shape investigation. 
Comput Aided Geom Des 13: 569--582.

\bibitem{ch}  Chern  S S (1945) 
Some new characterizations of the Euclidean sphere. 
Duke Math J 12:  279--290.

\bibitem{de}  Delaunay C (1841)
 Sur la surface de r\'evolution dont la courbure moyenne est
constante. J Math Pure Appl  6: 309--320.


\bibitem{en1}   Enneper A (1866)
 Ueber die cyclischen Fl\"{a}chen.
Nach K\"{o}nigl Ges  Wissensch G\"{o}ttingen Math Phys Kl: 243--249.

\bibitem{en2}  Enneper A (1869) 
Die cyklischen Fl\"{a}chen.  Z Math Phys 14:  393--421.

\bibitem{gmm}  G\'alvez J A, Mart\'{\i}nez  A, Milá\'an F (2003)
Linear Weingarten Surfaces in $\r^3$.   Monatsh Math 138: 133--144.

\bibitem{hw} Hartman P, Winter W (1954) 
 Umbilical points and W-surfaces. Am  J Math 76: 502-–508.

\bibitem{ho}  Hopf H (1951)
 \"{U}ber Fl\"{a}chen mit einer Relation zwischen den Hauptkr\"{u}mmungen. 
Math Nachr 4: 232--249.



\bibitem{ks}  K\"{u}hnel W,  Steller M (2005)
On closed Weingarten surfaces.  Monatsh Math 146: 113--126.


\bibitem{ni}   Nitsche J C C (1989)
 Cyclic surfaces of constant mean curvature.
Nachr Akad Wiss Gottingen Math Phys II 1: 1--5.

\bibitem{ri}  Riemann B (1868)
\"{U}ber die Fl\"{a}chen vom Kleinsten Inhalt be gegebener Begrenzung. 
Abh K\"{o}nigl Ges  Wissensch. G\"{o}ttingen Mathema. Cl.  13: 329--333.

\bibitem{rs} Rosenberg H,  Sa Earp R (1994)
The geometry of properly embedded special surfaces in $\r^3$; e.g., surfaces satisfying $aH+bK=1$, where $a$ and $b$ are positive.
 Duke Math. J. 73:  291-–306.


\bibitem{vo}  Voss K  (1959)
\"{U}ber geschlossene Weingartensche Fl\"{a}chen.
Math Annalen 138: 42--54.

 \bibitem{we1} Weingarten J (1861) 
Ueber eine Klasse auf einander abwickelbarer Fl\"{a}achen. 
J. Reine Angew. Math. 59:  382--393.

 \bibitem{we2}  Weingarten J (1863)   
Ueber die Fl\"{a}chen, derer Normalen eine gegebene Fl\"{a}che ber\"{u}hren.
 J. Reine Angew. Math.  62: 61-–63.



 
\end{thebibliography}
\end{document}